\documentclass[a4paper,11pt,oneside]{amsart}

\usepackage{hyperref}
\usepackage{amssymb,amsmath,amssymb,amsfonts,amsthm}
\usepackage{bm}
\usepackage{cases}
\usepackage{url}
\usepackage{geometry}
\usepackage{cite,mathtools}
\usepackage{xcolor}

\theoremstyle{plain}
\newtheorem{theorem}{Theorem}
\newtheorem{lemma}[theorem]{Lemma}
\newtheorem{problem}[theorem]{Problem}
\newtheorem*{algo}{\bf Algorithm}

\theoremstyle{definition}

\allowdisplaybreaks[4]

\def\BN{\mathbb N}
\def\BR{\mathbb R}
\def\T{\mathrm T}
\def\rd{\mathrm d}
\def\Ga{\Gamma}
\def\Om{\Omega}
\def\al{\alpha}
\def\be{\beta}
\def\de{\delta}
\def\ep{\epsilon}
\def\f{\frac}
\def\nb{\nabla}
\def\ov{\overline}
\def\pa{\partial}

\title[Numerical Reconstruction of Orders in Coupled Subdiffusion Systems]{Numerical Reconstruction of Orders in Coupled Systems of Subdiffusion Equations}

\author[Y. Liu]{Yikan Liu}
\address{Department of Mathematics, Kyoto University, Kitashirakawa-Oiwakecho, Sakyo-ku, Kyoto 606-8502, Japan.}
\email{liu.yikan.8z@kyoto-u.ac.jp}


\begin{document}

\maketitle

\begin{abstract}
In this paper, we develop a numerical algorithm for an inverse problem on determining fractional orders of time derivatives simultaneously in a coupled subdiffusion system. Following the theoretical uniqueness, we reformulate the order inverse problem as a discrete minimization problem, so that we derive a concise Gauss-Newton iterative method. Abundant numerical tests demonstrate the efficiency and accuracy of the proposed algorithm.
\end{abstract}

\section{Introduction}\label{sec-intro}

The last decade has witnessed explosive developments and gradual saturation of nonlocal models from various practical backgrounds. As a popular representative, partial differential equations with fractional derivatives in time have gathered considerable attentions among applied mathematicians. Recently, linear theory for subdiffusion equations with time derivatives in $(0,1)$ has been well established, followed by extensive researches on related numerical methods and inverse problems (e.g., \cite{JZ23,KR23,KRY20,SY11}). Nevertheless, it turns out that most existing literature is restricted to single equations, whose coupled counterparts have not been well studied yet.

Based on the above background, this article follows the line of \cite{LHL23} to focus on an inverse problem for coupled systems of subdiffusion equations. To formulate the problem, first we fix notations and define fractional derivatives. Throughout, let $T>0$ be constant and $\Om\subset\BR^d$ ($d\in\BN:=\{1,2,\dots\}$) be a bounded domain with a smooth boundary $\pa\Om$. For a constant $\al\in(0,1)$, let $\pa_t^\al$ be the algebraic inverse of the Riemann-Liouville integral operator $J^\al:L^2(0,T)\longrightarrow L^2(0,T)$ defined by
\[
J^\al f(t):=\int_0^t\f{(t-s)^{\al-1}}{\Ga(\al)}f(s)\,\rd s,\quad f\in L^2(0,T),
\]
where $\Ga(\,\cdot\,)$ denotes the Gamma function. Then according to \cite{KRY20}, the domain $D(\pa_t^\al)$ of $\pa_t^\al$ is $\ov{\{f\in C^1[0,T]\mid f(0)=0\}}^{H^\al(0,T)}$, where $H^\al(0,T)$ denotes the Sobolev-Slobodeckij space (see \cite{A75}). For a constant $K\in\BN$, let constants $\al_1,\dots,\al_K$ satisfy $1>\al_1\ge\cdots\ge\al_K>0$. In this paper, we consider the initial-boundary value problem
\begin{equation}\label{eq-ibvp}
\left\{\!\begin{alignedat}{2}
& \pa_t^{\al_k}\left(u_k-u_0^{(k)}\right)=\mathrm{div}(\bm A_k(\bm x)\nb u_k)+\sum_{\ell=1}^K c_{k\ell}(\bm x)u_\ell & \quad & \mbox{in }\Om\times(0,T),\\
& u_k=0 & \quad & \mbox{on }\pa\Om\times(0,T)
\end{alignedat}\right.
\end{equation}
for $k=1,\dots,K$, where $\bm A_k=(a_{ij}^{(k)})_{1\le i,j\le d}\in C^1(\ov\Om;\BR_{\mathrm{sym}}^{d\times d})$ ($k=1,\dots,K$) are strictly positive-definite matrix-valued functions on $\ov\Om$, and $c_{k\ell}\in L^\infty(\Om)$ ($k,\ell=1,\dots,K$). Here each $u_0^{(k)}$ stands for the initial value of the $k$-th component $u_k$, and the initial condition should be interpreted as $u_k(\bm x,\,\cdot\,)-u_0^{(k)}(\bm x)\in D(\pa_t^{\al_k})$ for a.e.\! $\bm x\in\Om$. Then the governing equations in \eqref{eq-ibvp} is a weakly coupled subdiffusion system with $K$ components, which are coupled with each other via coupling coefficients $c_{k\ell}(\bm x)$.

For later simplicity, we abbreviate $\bm\al=(\al_1,\dots,\al_K)^\T$ and $\bm u=(u_1,\dots,u_K)^\T$, where $(\,\cdot\,)^\T$ denotes the transpose. Now, we propose the order inverse problem under consideration in this article.

\begin{problem}\label{prob}
Let $\bm u$ satisfy \eqref{eq-ibvp} and fix $\bm x_0\in\Om$, $k_0\in\{1,\dots,K\}$ arbitrarily. Determine the orders $\bm\al$ in \eqref{eq-ibvp} by the single point observation of the $k_0$-th component $u_{k_0}$ of $\bm u$ at $\{\bm x_0\}\times(0,T)$.
\end{problem}

Due to practical importance, inverse problems on determining fractional orders in nonlocal equations have been investigated extensively from both theoretical and numerical aspects. We refer to \cite{LLY19} as a comprehensive survey before 2019, and recent developments on this topic involves the unique determination of orders with weaker data, e.g., in unknown media (see \cite{JK21,JK22}) or by short-time behavior (see \cite{LY23}). On the same direction, Problem \ref{prob} seeks the possibility of identifying all unknown orders simultaneously by the data of a single component by means of the coupling effect. As for other inverse problems for coupled subdiffusion systems, see \cite{RHY21} for a special coefficient inverse problem and \cite{FLL24} for a backward problem.

Based on the uniqueness for Problem \ref{prob} proved in \cite{LHL23} (see Lemma \ref{lem-ip}), this paper is devoted to the numerical reconstruction of all unknown orders in \eqref{eq-ibvp} by the same data, which seems not studied for coupled subdiffusion systems. In Section \ref{sec-recon}, we revisit the theoretical results and derive the numerical reconstruction method. Then, in Section \ref{sec-numer}, we implement the proposed method and demonstrate several numerical examples.


\section{Preliminary and Derivation of the Reconstruction Method}\label{sec-recon}

In this section, we first recall the well-posedness of the forward problem \eqref{eq-ibvp} and the uniqueness of Problem \ref{prob}. In \cite[Theorem 1]{LHL23}, the well-posedness of the mild solution was established in a more general formulation than that of \eqref{eq-ibvp} with $t$-dependent moderate coupling. Here, we summarize a simplified version in the following lemma.

\begin{lemma}\label{lem-fp}
Let $\bm u_0\in(L^2(\Om))^K$. Then there exists a unique mild solution $\bm u\in L^1(0,T;$ $(H^2(\Om)\cap H_0^1(\Om))^K)$ to \eqref{eq-ibvp}. Moreover, there exists a constant $C>0$ depending only on $\Om,\bm\al,T,\bm A_k,c_{k\ell}$ $(k,\ell=1,\dots,K)$ such that
\[
\|\bm u\|_{L^1(0,T;(H^2(\Om)\cap H_0^1(\Om))^K)}\le C\|\bm u_0\|_{(L^2(\Om))^K}.
\]
\end{lemma}

For the definition of a mild solution to \eqref{eq-ibvp}, see \cite{LHL23} for details. Recently, \cite{FLL24} shown the well-posedness of the strongly coupled version of \eqref{eq-ibvp} assuming the coincidence of all orders $\al_k$. As related studies on nonlinear generalizations of \eqref{eq-ibvp}, we refer to \cite{FY24,GW20,S21,S22}.

By the solution regularity in the above lemma and the Sobolev embedding $H^2(\Om)\subset C(\ov\Om)$ for the spatial dimensions $d=1,2,3$, we see that the observation data $u_{k_0}(\bm x_0,\,\cdot\,)$ makes sense in $L^1(0,T)$. Now we state the uniqueness result of Problem \ref{prob} obtained in \cite[Theorem 3]{LHL23}.

\begin{lemma}\label{lem-ip}
Let $d=1,2,3$ and $\bm u_0\in(L^2(\Om))^K$ satisfy $u^{(k)}_0\ge0,\not\equiv0$ in $\Om$ for all $k=1,\dots, K$. Assume that the coupling coefficients $c_{k\ell}$ satisfy
\begin{gather}
c_{k\ell}\ge0,\not\equiv0\quad\mbox{in }\Om,\ k,\ell=1,\dots,K,\ k\ne\ell,\label{eq-cond1}\\
\sum_{\ell=1}^K c_{k\ell}\le0\quad\mbox{in }\Om,\ k=1,\dots, K.\label{eq-cond2}
\end{gather}
Further, let $\bm u$ and $\bm v$ be the solutions to \eqref{eq-ibvp} with the fractional orders $\bm\al$ and $\bm\be,$ respectively. Then $u_{k_0}=v_{k_0}$ at $\{\bm x_0\}\times(0,T)$ implies $\bm\al=\bm\be$.
\end{lemma}

The non-negativity of the initial values and \eqref{eq-cond2} are technical assumptions, but \eqref{eq-cond1} means the cooperativeness and non-decoupling between components, which validates the simultaneous identification of all fractional orders by observing only a single component. Meanwhile, \eqref{eq-cond1}--\eqref{eq-cond2} implies that the decay rate $c_{kk}$ overwhelms the cooperativeness effect $\sum_{\ell\ne k}c_{k\ell}$ for each $k=1,\dots,K$, which results in the vanishing of all components as $t\to\infty$.

Based on the theoretical uniqueness stated above, we derive a numerical reconstruction method for Problem \ref{prob}. For simplicity, we restrict the spatial dimension $d=1$ and assume $A_k(x)\equiv1$, i.e., $\mathrm{div}(A_k(x)\nb u_k)=\pa_x^2u_k$ for all $k=1,\dots,K$. Denoting $\Om=(0,L)$ ($L>0$) and the coupling matrix $\bm C=(c_{k\ell})_{1\le k,\ell\le K}$, we deal with
\begin{equation}\label{eq-ibvp-sys3}
\begin{cases}
\pa^{\bm\al}_t(\bm u-\bm u_0)=\pa^2_x\bm u+\bm C\bm u & \mbox{in }(0,L)\times(0,T),\\
\bm u=\bm0 & \mbox{on }(0,L)\times(0,T).
\end{cases}
\end{equation}
We discretize the interval $(0,T)$ by an equidistant partition $0=t_0<t_1<\cdots<t_N=T$ with $t_i=\f TNi$ ($N\in\BN,i=0,1,\dots,N$). For $x_0\in(0,L)$, by $g_i$ we denote the observation data of $u_{k_0}(x_0,t_i)$. Then we simply reformulate Problem \ref{prob} as the following discrete minimization problem
\begin{equation}\label{eq-min}
\min_{\bm\al\in(0,1)^K}\Phi(\bm\al),\quad\Phi(\bm\al):=\f12\sum_{i=1}^N\left|u_{k_0}(\bm\al)(x_0,t_i)-g_i\right|^2.
\end{equation}
Here and henceforth, by $\bm u(\bm\al)=(u_1(\bm\al),\dots,u_K(\bm\al))^\T$ we emphasize the dependency of the solution to \eqref{eq-ibvp-sys3} on the order $\bm\al$. Further, introducing the residuals
\[
r_i(\bm\al)=u_{k_0}(\bm\al)(x_0,t_i)-g_i\ (i=1,\dots,N),\quad\bm r(\bm\al)=(r_1(\bm\al),\dots,r_N(\bm\al))^\T,
\]
we can represent the Jacobian matrix of $\Phi(\bm\al)$ as
\[
\bm J(\bm\al)=
\begin{pmatrix}
\f{\pa r_1}{\pa\al_1}(\bm\al) &\cdots &\f{\pa r_1}{\pa\al_K}(\bm\al)\\
\vdots & \ddots & \vdots\\
\f{\pa r_N}{\pa\al_1}(\bm\al) &\cdots &\f{\pa r_N}{\pa\al_K}(\bm\al)
\end{pmatrix}.
\]
Therefore, it follows from direct calculations that the gradient of $\Phi(\bm\al)$ can be written as $\nb\Phi(\bm\al)=\bm J(\bm\al)^\T\bm r(\bm\al)$.

Until now, we have discussed the formulation of Problem \ref{prob}, that is, observing only a single component. For the numerical reconstruction, we are also interested in comparing the numerical performance of observing multiple components, for which the same argument as above also works. For instance, in the case of observing 3 components $k_0,k_1,k_2\in\{1,\dots,K\}$, the discrete minimization problem corresponding to \eqref{eq-min} turns out to be
\[
\min_{\bm\al\in(0,1)^K}\Phi(\bm\al),\quad\Phi(\bm\al):=\f12\sum_{j=0}^2\sum_{i=1}^N\left|u_{k_j}(\bm\al)(x_0,t_i)-g_i^{(k_j)}\right|^2,
\]
where $g_i^{(k_j)}$ stands for the observation data of $u_{k_j}(x_0,t_i)$ ($j=0,1,2$, $i=1,\dots,N$). Further, setting the residuals as
\begin{gather*}
r_i^{(k_j)}(\bm\al)=u_{k_j}(\bm\al)(x_0,t_i)-g_i^{(k_j)}\quad(j=0,1,2,\ i=1,\dots,N),\\
\bm r(\bm\al)=\left(r_1^{(k_0)}(\bm\al),\dots,r_N^{(k_0)}(\bm\al),r_1^{(k_1)}(\bm\al),\dots,r_N^{(k_1)}(\bm\al),r_1^{(k_2)}(\bm\al),\dots,r_N^{(k_2)}(\bm\al)\right)^\T,
\end{gather*}
we can readily verify that the expression $\nb\Phi(\bm\al)=\bm J(\bm\al)^\T\bm r(\bm\al)$ keeps the same.

In summary, we propose the following iterative algorithm for Problem \ref{prob} taking advantage of the Gauss-Newton method.

\begin{algo}
Fix the tolerance $\ep>0$ for convergence, the initial guess $\bm\al_0\in(0,1)^K$ and set $m=0$.

{\bf Step 1} Denote $\bm J_m=\bm J(\bm\al_m)$ and compute $\bm s_m=-(\bm J^\T_m \bm J_m)^{-1}\bm J^\T_m\bm r(\bm\al_m)$.

{\bf Step 2} Update $\bm\al$ by $\bm\al_{m+1}=\bm\al_m+\bm s_m$.

{\bf Step 3} If $|\bm s_m|\le\ep$, then return $\bm\al_*=\bm\al_{m+1}$ as the reconstructed solution and stop the iteration. Otherwise, update $m\leftarrow m+1$ and return to Step 1.
\end{algo}

In each step of the iteration, obviously the major computational cost lies in the approximation of the Jacobian matrix $\bm J_m$. To compute the component $\f{\pa r_i}{\pa\al_k}(\bm\al)$ in $\bm J(\bm\al)$, we add a small perturbation to the $k$-th component of $\bm\al$, i.e.,
\[
\bm\al+\ep\bm e_k=(\al_1,\dots,\al_{k-1},\al_k+\ep,\al_{k+1},\dots,\al_K)^\T
\]
with some $|\ep|\ll1$ and approximate $\f{\pa r_i}{\pa\al_k}(\bm\al)$ by the difference quotient
\[
\f{\pa r_i}{\pa\al_k}(\bm\al)\approx\f{u_{k_0}(\bm\al+\ep\bm e_k)-u_{k_0}(\bm\al)}\ep.
\]
As a result, in each step we should solve $K+1$ forward problems \eqref{eq-ibvp-sys3} with slightly perturbed orders $\bm\al$. As for the numerical simulation of \eqref{eq-ibvp-sys3}, the fractional derivatives in time is discretized by the L1 scheme proposed in \cite{LX07}, and the spatial derivative $\pa_x^2$ is discretized by the traditional central difference quotient. Then it is not difficult to deduce a numerical scheme in form of linear systems, and we omit the details here. The forward solver turns out to be efficient enough, so that the computational cost of the above proposed iterative algorithm is definitely affordable at least for $d=1$.


\section{Numerical Examples}\label{sec-numer}

In this section, we test the performance of the above algorithm via abundant numerical examples.

We fix the basic settings of the numerical reconstruction. In the initial-boundary value problem \eqref{eq-ibvp-sys3}, we take $L=T=1$ and $N=100$ for the discretization in time. We choose $x_0=0.5$ as the observation point, and the noisy observation data $g_i^\de$ are constructed randomly with the noise level $\de>0$ such that $|g_i^\de-g_i|\le\de g_i$ ($i=1,\dots,N$).

We start with the simplest coupled case of \eqref{eq-ibvp-sys3}, that is, $K=2$. We set the true solution as $\bm\al=(0.9,0.5)^\T$, and the coupling coefficients and initial values as
\[
\begin{gathered}
c_{11}=c_{22}=-1,\quad c_{12}=c_{21}=1,\\
u_0^{(1)}(x)=2u_0^{(2)}(x),\quad u_0^{(2)}(x)=-4\left(x-\f12\right)^2+1\ge0,
\end{gathered}
\]
respectively. For the components of observation, we study all 3 possible combinations, that is,

{\bf Case (A)} $u_1(x_0,t)$,

{\bf Case (B)} $u_2(x_0,t)$,

{\bf Case (C)} $u_1(x_0,t)$ and $u_2(x_0,t)$.

We first test the algorithm with the initial guess $\bm\al_0=(0.5,0.5)^\T$, whose second component coincides with that of the true solution. The numerical results and the relative errors with 3 choices of the noise level $\de=0\%,1\%,5\%$ are listed in Tables \ref{table-num-cal-results-0505} and \ref{table-relative-errors-0505}, respectively. Especially in Table \ref{table-relative-errors-0505}, we see that in all cases the related errors are lower than $0.05\%$ for $\de=0\%$ (noiseless data) and lower than $0.7\%$ for $\de=5\%$. Further, the errors of observing $u_2$ (Case (B)) are slightly larger than other cases for noisy data, while there seems no essential difference between observing $u_1$ (Case (A)) and both $u_1,u_2$ (Case (C)). All these results demonstrate the high accuracy of the proposed numerical method even without employing any advanced regularization techniques. The reason partially lies in the fact that Problem \ref{prob} is highly overdetermined, namely, we use at least one $t$-dependent function to find only 2 constants. Therefore, although observing both components doubles the amount of data, it does not greatly improve the numerical performance. On the other hand, all iterations terminate within 5 steps, which inherits the rapid convergence of Newton-type iterative methods.
\begin{table}[htbp]\centering
\caption{The numerical results of reconstructing $\bm\al$ with the initial guess $\bm\al_0=(0.5, 0.5)^\T$.}\label{table-num-cal-results-0505}
\begin{tabular}{c|ccc}
\hline\hline
$\de$ & 0\% & 1\% & 5\%\\
\hline
Case (A) & (0.899724, 0.499772) & (0.899274, 0.500281) & (0.897462, 0.502305)\\
Case (B) & (0.899723, 0.499823) & (0.899156, 0.500521) & (0.896872, 0.503315)\\
Case (C) & (0.899724, 0.499772) & (0.899273, 0.500283) & (0.897460, 0.502318)\\
\hline\hline
\end{tabular}\bigskip

\caption{The relative errors ($\%$) of reconstructing $\bm\al$ with the initial guess $\bm\al_0=(0.5, 0.5)^\T$.}\label{table-relative-errors-0505}
\begin{tabular}{c|ccc}
\hline\hline
$\de$ & 0\% & 1\% & 5\%\\
\hline
Case (A) & (0.030665, 0.045580) & (0.080719, 0.056110) & (0.281960, 0.460937)\\
Case (B) & (0.030756, 0.035474) & (0.093815, 0.104164) & (0.347606, 0.663082)\\
Case (C) & (0.030665, 0.045574) & (0.080771,0.056632) & (0.282228, 0.463542)\\
\hline\hline
\end{tabular}
\end{table}

Since the initial guess $\bm\al_0=(0.5,0.5)^\T$ partly coincides with the true solution, we test the the algorithm with $\bm\al_0=(0.7,0.3)^\T$ to confirm its effectiveness. We repeat the same procedure and list the relative errors in Table \ref{table-relative-errors-0703}. As expected, again we see the satisfactory performance of the method.
\begin{table}[htbp]\centering
\caption{The relative errors ($\%$) of reconstructing $\bm\al$ with the initial guess $\bm\al_0=(0.7, 0.3)^\T$.}\label{table-relative-errors-0703}
\begin{tabular}{c|ccc}
\hline\hline
$\de$ & 0\% & 1\% & 5\%\\
\hline
Case (A) & (0.030768, 0.034912) & (0.080822, 0.065833) & (0.282058, 0.467454)\\
Case (B) & (0.030777, 0.033635) & (0.093837, 0.105731) & (0.347614, 0.662978)\\
Case (C) & (0.030768, 0.034911) & (0.080875, 0.066348) & (0.282326, 0.470047)\\
\hline\hline
\end{tabular}
\end{table}

Since Problem \ref{prob} is nonlinear, the minimization problem \eqref{eq-min} is not convex and thus the algorithm may not converge. To this end, we test the convergence of the method with initial guesses $\bm\al_0=(i/10,j/10)$ for $9\ge i\ge j\ge1$ ($i,j\in\BN$). Then there are 45 choices of $(i,j)$ and we run the algorithm for all Cases (A)--(C). In Table \ref{table-conv1}, we record the times of convergence / divergence and the iteration steps if it converges. Again we see the convergence to the true solution within 5 steps in all cases, indicating the high efficiency of the algorithm. Meanwhile, the convergence occurs in major choices of initial guesses. Indeed, we observe that if the first component of $\bm\al_0$ is far away from the true solution $0.9$, say, between $0.1$ and $0.4$, the algorithm mostly diverges.
\begin{table}[htbp]\centering
\caption{Situations of convergence / divergence in 45 choices of initial guesses ($K=2$).}\label{table-conv1}
\begin{tabular}{c|ccc}
\hline\hline
& convergence & divergence & iteration steps\\
\hline
Case (A) & 39 & 6 & $1\sim5$\\
Case (B) & 35 & 10 & $1\sim5$\\
Case (C) & 39 & 6 & $1\sim5$\\
\hline\hline
\end{tabular}
\end{table}

Finally, we proceed to the case of $K=3$ and set
\[
\begin{gathered}
c_{kk}=-2\ (k=1,2,3),\quad c_{k\ell}=1\ (k,\ell=1,2,3,\ k\ne\ell),\\
u_0^{(1)}(x)=2u_0^{(2)}(x),\quad u_0^{(2)}(x)=u_0^{(3)}(x)=-4\left(x-\f12\right)^2+1\ge0.
\end{gathered}
\]
We choose the true solution as $\bm\al=(0.9,0.6,0.5)^\T$ and observe the following 3 out of 6 combinations of components:

{\bf Case (D)} $u_3(x_0,t)$,

{\bf Case (E)} $u_2(x_0,t)$ and $u_3(x_0,t)$,

{\bf Case (F)} $u_1(x_0,t)$, $u_2(x_0,t)$ and $u_3(x_0,t)$.

\noindent As before, we change the initial guesses as $\bm\al_0=(i/10,j/10,k/10)$ for $9\ge i\ge j\ge k\ge1$ ($i,j,k\in\BN$), so that there are altogether 165 choices. We test the algorithm in all cases and record the times of convergence / divergence and the iteration steps in Table \ref{table-conv2}. Again we verify the efficiency and accuracy of the proposed method, which, however, yields more divergent results than that for $K=2$. Indeed, again we observe that the algorithm diverges if the first component of $\bm\al_0$ is far away from the true solution $0.9$, say, between $0.1$ and $0.7$.
\begin{table}[htbp]\centering
\caption{Situations of convergence / divergence in 165 choices of initial guesses ($K=3$).}\label{table-conv2}
\begin{tabular}{c|ccc}
\hline\hline
& convergence & divergence & iteration steps\\
\hline
Case (D) & 29 & 136 & $1\sim4$\\
Case (E) & 78 & 87 & $1\sim6$\\
Case (F) & 77 & 88 & $1\sim5$\\
\hline\hline
\end{tabular}
\end{table}

The numerical tests for components $K>3$ and / or spatial dimensions $d>1$ require higher computational costs, which await future investigations.\bigskip

{\bf Acknowledgements } This work is supported by JSPS KAKENHI Grant Numbers JP22K13954 and JP23KK0049.


\end{document}